\documentstyle{amsart}
\input{bull-art}
\newcommand{\R}{{\bold R}}
\newcommand{\Rn}{{\bold R}^n}
\newcommand{\Z}{{\bold Z}}
\newcommand{\Sc}{{\scr S}}

\newcommand{\Pj}{P_j}

\renewcommand{\l}{\lambda}

\newcommand{\e}{\epsilon}
\newcommand{\ov}{\overline}

\newtheorem{theorem}{Theorem}[section]
\newtheorem{proposition}[theorem]{Proposition}

\newtheorem{corollary}[theorem]{Corollary}

\theoremstyle{definition}
\newtheorem{definition}[theorem]{Definition}

\newcommand{\abs}[1]{\left|{#1}\right|}

\newcommand{\F}{\phi}
\newcommand{\M}{\psi}

\newcommand{\Fjk}{\phi_{jk}}

\newcommand{\Mjk}{\psi_{jk}}
\newcommand{\Fx}{\phi (2^j x-k)}

\newcommand{\RB}{ {\scr R}{\scr B} }

\begin{document}
\def\currentvolume{30}
\def\currentissue{1}
\def\currentyear{1994}
\def\currentmonth{January}
\def\copyrightyear{1994}
\def\currentpages{87-94}

\title{Pointwise Convergence of Wavelet Expansions}
\ratitle
    \author[S. E. Kelly]{Susan E. Kelly}

\thanks{The first and second authors' research was 
partially supported by  
the National Science Foundation} 
\thanks{The second and third authors' research was
partially supported by the Air Force Office of Scientific 
Research} 
\author[M. A. Kon]{Mark A. Kon} 
\author[L. A. Raphael]{Louise A. Raphael}
\address[S. Kelly]{Department of Mathematics, University 
of Wisconsin,
La Crosse, La Crosse, WI 54601}
\email{kelly@@math.uwlax.edu}
\address[M. Kon]{Department of Mathematics, Boston 
University,
Boston, MA 02215}
\email{mkon@@math.bu.edu}
\address[L. Raphael]{Department of Mathematics, Howard 
University,
Washington, DC 20059}
\subjclass{Primary 42C15; Secondary 40A30}
\date{July 9, 1993}
\thanks{Some contents of this article were presented in 
June 1992
at the Joint Summer Research Conference on Wavelets and 
Applications,
Mt.\ Holyoke College, South Hadley, Massachusetts}

\maketitle

\begin{abstract}
In this note we announce that under general hypotheses, 
wavelet-type
expansions (of functions in $L^p,\ 1\leq p \leq \infty$, 
in one or
more dimensions) converge pointwise almost everywhere, and 
identify
the Lebesgue set of a function as a set of full measure on 
which they
converge.  It is shown that unlike the Fourier summation 
kernel,
wavelet summation kernels $P_j$ are bounded by radial 
decreasing $L^1$
convolution kernels.  As a corollary it follows that best 
$L^2$
spline approximations on uniform meshes converge
pointwise almost everywhere.  Moreover, summation of wavelet
expansions is partially insensitive to order of summation.

We also give necessary and sufficient conditions for given 
rates of
convergence of wavelet expansions in the sup norm.  Such 
expansions
have order of convergence $s$ if and only if the basic 
wavelet $\psi$
is in the homogeneous Sobolev space $H^{-s-d/2}_h$.  We 
also present
equivalent necessary and sufficient conditions on the 
scaling
function.  The above results hold in one and in multiple 
dimensions.
\end{abstract}


\section{Introduction}

In this note we present several convergence results for 
wavelet and
multi\-resol\-ution-type expansions.  It is natural to ask 
where such
expansions converge (and whether they converge almost 
everywhere) and
what are the rates of convergence.  The answer is that 
under rather
weak hypotheses, one- and multidimensional 
wavelet expansions
converge pointwise almost everywhere and, more 
specifically, on the
Lebesgue set of a function being expanded.  We will also 
give exact
rates of convergence in the supremum norm, in terms of 
Sobolev
properties of the basic wavelet or of the scaling function.

Wavelets with local support in the time and frequency 
domains were
defined by A.~Grossman and J.~Morlet \cite{Gr-Mo1} in 1984 
in order
to analyze seismic data.  The prototypes of wavelets, 
however, can be
found in the work of A.~Haar \cite{Ha1} and the modified 
Franklin
systems of J.-O.~Str\"omberg \cite{St1}.

To identify the underlying structure and to generate
interesting examples of orthonormal bases for $L^2(\bold 
R)$, S.~Mallat
\cite{Ma1} and Y.~Meyer \cite{Me1} developed the framework 
of multiresolution
analysis.  P.~G. Lemari\'{e} and Y.~Meyer \cite{Le-Me1} 
constructed
wavelets in $\Sc (\Rn )$, the space of rapidly decreasing 
smooth functions.
J.-O. Str\"{o}mberg \cite{St1} developed spline wavelets 
while looking
for unconditional bases for Hardy spaces.  G.~Battle 
\cite{Ba1} and
P.~G. Lemari\'{e} \cite{Le1} developed these bases in the 
context of
wavelet theory.  Spline wavelets have exponential decay 
but only
$C^N$ smoothness (for a finite $N$ depending on the order 
of the
associated splines).  I.~Daubechies \cite{Da1} constructed 
compactly
supported wavelets with $C^N$ smoothness.  The support of 
these
wavelets increased with the smoothness; in general, to 
have $C^\infty
$ smoothness, wavelets must have infinite support.

Y.~Meyer \cite{Me1} was among the first to study 
convergence results
for wavelet expansions; he was followed by G.~Walter 
\cite{Wa1,
Wa2}.  Meyer proved that under some regularity assumptions
on the wavelets, wavelet expansions of continuous 
functions converge
everywhere. In contrast to these results, the pointwise 
convergence
results presented here give almost everywhere convergence 
(and
convergence on the Lebesgue set) for expansions of general 
$L^p$
($1\leq p\leq \infty$) functions.  We assume rather mild 
bounds and no
differentiability for the wavelet or the scaling function; 
our
conditions allow inclusion of the families of so-called 
$r$-regular
wavelets \cite{Me1}, as well as some other wavelets.

These results parallel L.~Carleson's \cite{Ca1} and 
R.~A.~Hunt's
\cite{Hu1} theorems for Fourier series.  One difference 
and slight
advantage of wavelet expansions comes from the fact that 
almost
everywhere convergence occurs on a simple set of full 
measure, namely
the Lebesgue set, while almost everywhere convergence for 
Fourier
series is established on a much more elaborate set of full 
measure.  

We also give necessary and sufficient conditions for given 
pointwise
(sup-norm) rates of convergence of wavelet or 
multiresolution
expansions, in terms of Sobolev conditions on the basic 
wavelet or the
scaling function.  It has been shown previously by Mallat 
\cite{Ma1}
and Meyer \cite{Me1} that the Sobolev class of a function is
determined by the $L^2$ rates of convergence of its 
wavelet expansion.
Necessary and sufficient conditions for $L^2$ rates of 
convergence
which are analogous to our sup-norm conditions have been 
obtained by
de Boor, DeVore, and Ron \cite{Bo-De-Ro}, who have also 
studied
sup-norm convergence in more general situations
\cite{Bo-Ro}.

Our results on convergence rates can be viewed as a 
sharpening in the
context of wavelets of the well-known Strang-Fix 
\cite{St-Fi}
conditions for convergence of multiscale expansions.

The results given here hold for multiresolution expansions 
in
multiple dimension.  The proofs, which will appear elsewhere
\cite{Ke-Ko-Ra1}, involve the kernels of the partial
sums of such expansions and the above-mentioned result 
that such
kernels are bounded by rescalings of $L^1$ radial 
functions.  We
should add that such bounds for wavelet expansions are 
nontrivial and
arise from cancellations which occur in the sum 
representations of the
partial sum kernels.  Naive bounding of the summation 
kernels by using
absolute values in their sum representations fails to 
yield the needed
radial bounds for any class of wavelets.  We remark that 
such bounds
on the summation kernel can be obtained more easily by 
writing it
using the orthonormal translates $\phi(x-k)$ of the 
scaling function,
instead of the wavelets $\psi_{jk}$.  However, in proving 
results for
convergence of wavelet expansions (Theorem 2.1(iii)), we 
wish to
avoid making any assumptions about radial bounds for the 
scaling
function.

The pointwise and $L^p$ convergence results contained here 
were
obtained independently by the first author and by the 
second two
authors.  Results on the Gibbs effect obtained by the 
first author
will appear elsewhere.

To start, we define a multiresolution analysis on 
$L^2(\R^d)$
\cite{Ma1, Me1}.
\begin{definition}
A {\em multiresolution analysis} on $L^2(\bold R^d)$, 
$d\geq 1$, is an
increasing
sequence $\{ V_j\}_{j\in\Z }$,
$$
\cdots\subset V_{-2}\subset V_{-1}\subset V_0\subset 
V_1\subset  
V_2\subset\cdots
$$
of closed subspaces in $L^2(\R^d )$ where
$$
\bigcap_{j\in \Z} V_j =\{ 0\}\, ,\qquad \ov{\bigcup_{j\in 
\Z} V_j} =  
L^2(\R^d )\, ,
$$
and the spaces $V_j$ satisfy the following additional 
properties:

\noindent (i) For all $f \in L^2(\R^d)$, $j\in \Z$, and 
$k\in \Z^d$,
$$
f(x) \in V_{j} \Longleftrightarrow f(2 x ) \in V_{j+1}\, 
$$
and
$$
f(x)\in V_0 \Longleftrightarrow f(x-k)\in V_0 \, .
$$

\noindent(ii) There exists a {\em scaling function} $\F 
\in V_0$
such that $\{\Fjk \}_{k\in\Z^d}$ is an orthonormal 
basis of $V_j$, where
\begin{equation}
\Fjk (x) = 2^{jd/2} \Fx\, ,
\end{equation}
for $x \in\R^d$, $j\in\Z$, and $k\in\Z^d$.

Associated with the $V_j$ spaces, we additionally define 
$W_j$ to be the
orthogonal complement of $V_j$ in $V_{j+1}$, so that
$V_{j+1} =  V_j \oplus W_j$.  Thus, $L^2(\R^d ) = 
\overline{\sum \oplus W_j}$.  
We define $P_j$ and $Q_j = P_{j+1}-P_j$, respectively, to 
be the  
orthogonal projections onto the spaces $V_j$ and $W_j$, 
with kernels  
$P_j(x,y)$ and $Q_j(x,y)$.
\end{definition}

Under the assumptions in the above definition and with some
additional regularity, it can be proved
\cite{Me1, Da2} that there then exists a
set $ \{\M^{\l}\} \in W_0$, where $\l$ belongs to an index 
set
$\Lambda $ of cardinality $2^d -1$, such that $\{ \Mjk^{\l}
\}_{k\in \Z^d,\lambda}$ is an orthonormal basis of $W_j$, 
and
thus $\{\Mjk^{\l} \}_{j\in\Z,k\in\Z^d,\lambda}$ is a {\it 
wavelet basis} of
$L^2(\R^d )$, where \begin{equation} \Mjk^{\l}(x) = 
2^{jd/2} \M^{\l}(2^jx-k)\, ,
\end{equation} for $x\in\R^d$, $j \in \Z$, and $k\in\Z^d$.

\begin{definition} 

For $f\in L^p(\R^d )$ ($1\leq p\leq\infty$) we define the 
following
related expansions:

(a) The sequence of projections $\{\Pj f(x)\}_j$ will be
called the {\em multiresolution expansion} of $f$.

(b) The {\em scaling expansion} of $f$ is defined to be
\begin{equation}
f\sim\sum_k b_k\F (x-k) + \sum_{j\geq 0;k;\l} 
a_{jk}\M^{\l}_{jk}(x) \, ,  
\label{BK}
\end{equation}
where the coefficients $a_{jk}$ and $b_k$ are the $L^2$ 
expansion  
coefficients of $f$. 

(c) The {\em wavelet expansion} of $f$ is
\begin{equation}
f\sim \sum_{j;k;\l} a_{jk}\M^{\l}_{jk} (x) \, ,    
\label{AJK}
\end{equation}
where the coefficients $a_{jk}$ are the $L^2$ expansion  
coefficients of $f$.
\end{definition}

We remark for part (a) of the above definition that it can 
be shown
that the projections $P_j$ (defined by their integral 
kernels) extend
to bounded operators on $L^p$, $1\leq p \leq \infty$.  The 
$L^2$
expansion coefficients in (b) and (c) (defined by 
integration
against $f$) are defined and uniformly bounded for any 
$f\in L^p$,
$1\leq p \leq \infty$.

\begin{definition} The point $x$ is a {\em Lebesgue
point} of the measurable function $f(x)$ on ${\bold R}^d$ 
if $f$ is
integrable in some neighborhood of $x$, and
$$\lim_{\epsilon\rightarrow 0}{1\over {{V(B_\epsilon  
)}}}\int_{B_\epsilon} |f(x)-f(x+y)|\,dy =0,$$
where $B_\epsilon$ denotes the ball of radius $\epsilon$ 
about the
origin, and $V$ denotes volume.
\end{definition}

Such points $p$ are essentially characterized by the fact 
that the
average values of $f(x)$ around $p$ converge to the values 
of
$f$ at these points, as averages are taken over smaller 
balls centered
at $x$.  Note that all continuity points are Lebesgue 
points, but
the converse is not true.

\begin{definition} 

A function $f(x)$ is in the class $\RB$ if it is
absolutely bounded by an $L^1$ radial decreasing function 
$\eta (x)$,
i.e., $\eta (x_1)=\eta (x_2)$ whenever $|x_1|=|x_2|$, 
$\eta (x_1)\leq\eta (x_2)$ when $|x_1|\geq |x_2|$, and $\eta
(x)\in L^1({\bold R}^d)$.

A function $f$ is {\em partially continuous} if there 
exists a set $A$ of  
vectors $a\in\R^d$ with positive measure such that  
$\lim_{\e\rightarrow 0}\M (x+\e a) = \M (x)$ for $a\in A$.
\end{definition}

\begin{definition}
The {\em homogeneous Sobolev space} of order $s$ is 
defined by
\begin{equation}
H_h^s\equiv\left\{ f\in L^2({\bold R}^d)\colon\ \Vert
f\Vert_{h,s}\equiv \sqrt{\int |\hat f(\xi )|^2 |\xi 
|^{2s}\,d\xi}
<\infty \right\}.\label{hhs}
\end{equation}
\end{definition}

The ordinary Sobolev space $H^s$ is defined as $H^s_h$ in 
(\ref{hhs}),
with replacement of $|\xi|^{2s}$ by $(|\xi|^{2s} + 1)$.  
Under the
Fourier transform, the space $H_h^s$ is a dense subspace 
of the
complete weighted $L^2$ space of {\it all} measurable 
$\hat f(\xi)$
with $\Vert f\Vert_{h,s}<\infty$.  This dense subspace 
consists of
those functions $\hat f(\xi)$ which are also in the 
regular unweighted
$L^2$ space.

Convergence rates of wavelet expansions are sensitive to 
both the
smoothness of the wavelet and the smoothness of the 
function being
expanded.  For a wavelet $\psi$ of given smoothness, the 
sensitivity
to the Sobolev space of the function being expanded 
disappears when
the function's Sobolev parameter is sufficiently large.

\begin{definition}
A family $\psi^\lambda$ of wavelets yields {\em pointwise 
order of
approximation} (or {\em pointwise order of convergence}) 
$r$ in the
space $H^s$ if for 
any function $f\in H^s$, the $j\text{th}$ order wavelet 
approximation
$P_jf$ satisfies
\begin{equation}
\Vert P_jf-f\Vert_\infty =O(2^{-jr}), \label{PMF}
\end{equation}
as $j$ tends to infinity.
More generally, the wavelets $\psi^\lambda$ yield 
pointwise order of
approximation (or convergence) $r$ if for any function $f$ 
which is
sufficiently smooth (i.e., is in a Sobolev space of 
sufficiently
large order $s$), (\ref{PMF}) holds.

\end{definition}

We will give several necessary and sufficient conditions (in
terms of their Fourier transforms and membership in 
homogeneous
Sobolev spaces) on the basic wavelet $\psi$ or the scaling 
function
$\phi$ for given orders of convergence. In practice, 
sufficient
smoothness for a function $f$ (in the sense of the above 
definition) will mean
that $f$ is in $H^{s+d/2}$ or a higher Sobolev space.


\section{Pointwise convergence results}  \label{results}

With the background given, our main results can now be 
presented.
\begin{theorem}                            \label{main}
\rom{(i)} Assume only that the scaling
function $\phi$ of a given multiresolution analysis is in 
$\RB$, i.e.,
that it is bounded by an $L^1$ radial decreasing function. 
Then for  
any
$f\in L^p({\bold R}^d)$ $(1\leq p\leq \infty)$, its 
multiresolution
expansion $\{P_jf\}$ converges to $f$ pointwise almost 
everywhere.

\rom{(ii)} If $\phi,\psi^\lambda\in \RB$ for all 
$\lambda$, then
also both the scaling
\rom{(}\ref{BK}\rom{) (}if $1\leq p\leq \infty)$ and 
wavelet 
\rom{(}\ref{AJK}\rom{) 
(}if $1\leq p < \infty)$ expansions of any $f\in L^p$ 
converge to $f$
pointwise almost everywhere.  If further $\psi^\lambda$ 
and $\phi$ are
\rom{(}partially\/\rom{)} continuous, then both expansions 
converge to $f$ on its
Lebesgue set.

\rom{(iii)} If we assume only $\psi^\lambda (x)\ln (2+
|x|)\in  
\RB$ for all $\lambda$, then for $f\in L^p$, its wavelet
\rom{(}for $1\leq p < \infty)$ and multiresolution 
\rom{(}for $1 \leq p \leq
\infty$\rom{)} expansions converge to $f$ pointwise almost 
everywhere.  If
further $\psi^\lambda$ is \rom{(}partially\/\rom{)} 
continuous for all $\lambda$,
then both the wavelet and multiresolution expansions 
converge to $f$
on its Lebesgue set.

\rom{(iv)} The last two statements hold for any order of 
summation
in which the range of the values of $j$ for which the sum 
over $k$ and
$\lambda$ is partially complete always remains bounded.

\end{theorem}

In statement (iv) above, the summation over $k$ and 
$\lambda$ is
{\it partially complete} for a fixed $j$ if it contains 
some terms,
but not all with the given value of $j$.  By the {\it 
range of values}
for which the sum is partially complete we mean the 
difference of the
largest and smallest value of $j$ for which the sum is 
partially
complete.  Statement (iv) requires that this range always 
be smaller
than some constant $M$.

The above result on convergence of multiresolution 
expansions applies
to spline expansions as well.  Given a uniform grid $K$ in 
$\R$, one
might ask whether given a function $f\in L^2(\R^d)$, the 
best $L^2$
approximations $P_j f$ of $f$ (by splines of a fixed 
polynomial order
$k$) converge to $f$ pointwise as the grid size goes to 0. 
 The answer
to this is affirmative.

\begin{corollary} 
For $f\in L^p(\R)$ $(1\leq p \leq \infty)$, the order $k$ 
best $L^2$
spline approximations $P_jf$ of $f$ converge to $f$ 
pointwise almost
everywhere, and more specifically on the Lebesgue set of 
$f$, as the
uniform mesh size goes to \rom{0}.

\end{corollary}

The proof follows from the fact that best spline 
approximations are
partial sums of multiresolution expansions, with some 
radially bounded
($\RB$) scaling function $\phi$.  This result also extends 
to
multidimensional splines.  Technically, the best $L^2$ 
approximation
of $f\in L^p$ only makes sense when $p = 2$, but it can be 
defined for functions in $L^p$ by continuous extension of 
the
projections $P_j$ from $L^2$ to $L^p$.

The following proposition is a consequence of the proof of 
Theorem
\ref{main}.  It has been proved before under somewhat 
stronger
hypotheses, yielding stronger conclusions in \cite{Me1}.

\begin{proposition}
Under the hypotheses of case \rom{(i)}, case \rom{(ii)}, 
or case 
\rom{(iii)} of Theorem
\ref{main}, $L^p$ convergence of the expansions also 
follows for
$1\leq p<\infty$.
\end{proposition}

Thus for wavelet series and more generally for one- and
multidimensional multiresolution expansions, essentially 
all hoped for
convergence properties hold, regardless of rates of 
convergence.

The basis for Theorem \ref{main} is the bound on the 
kernel of the  
projection
$P_j$ onto the scaling space $V_j$.  It can be shown that 
under any
of the hypotheses in Theorem 2.1, the kernel
$P_m(x,y)$ has the form
$$
P_m(x,y)=\sum_{j<m;k\in\Z^d;\lambda}  
\psi_{jk}^\lambda (x)\overline{\psi_{jk}^\lambda
(y)} = \sum_{k\in\Z}\F_{mk}(x)\overline{\F_{mk}}(y)\,
$$
for $x, y \in {\bold R}^d$, with convergence of both sums 
on the right
occurring pointwise, uniformly on subsets a positive 
distance away
from the diagonal $D=\{(x,y):\ x=y\}$.
The kernel converges to a delta distribution $\delta 
(x-y)$ in the
following sense:

\begin{theorem}
Under the assumption that $\phi\in \RB$ 
or that $\psi^{\l} (x)\ln (2+|x|)\in \RB$ for all 
$\lambda$, the
kernels $P_j(x,y)$ of the projections onto $V_j$ satisfy the
convolution bound
\begin{equation}
|P_j(x,y)|\leq C2^{jd}H(2^{j}(|x-y|)), \label{pmx}
\end{equation}
where $H(|\cdot|)\in \RB$, i.e., $H(|\cdot|)$ is a radial 
decreasing
$L^1$ function. 
\end{theorem}

In one dimension, precise bounds can be obtained for kernels
of specific wavelets. Two examples are illustrated in the 
result
below.

\begin{theorem}
In $\R^1$, let $P_j(x,y)=\sum_{k\in {\bold 
Z}}\phi_{jk}(x)\phi_{jk}(y)$
be the summation kernel, generated by the scaling function 
$\phi \in
L^2(\R)$.

\rom{(a)} If $\phi$ has exponential decay, i.e., $\phi (x) 
\leq 
Ce^{-a|x|}$ for some positive $a$, then $$ |P_j(x,y)|\leq C
2^je^{-a2^j|x-y|/2}.  
$$

\rom{(b)} If $\phi$ has algebraic decay, $\phi (x)\leq 
{{C_N}\over
{{(1+|x|)^N}}}$ for some $N>1$, then 
$$
|P_j(x,y)|\leq C_N {{2^j}\over {{(1+2^j|x-y|)^N}}}\leq C_N 
2^j
$$ 
for $N>1$.

\end{theorem}

As a corollary, if a scaling function $\phi$ has rapid 
decay (faster
than any polynomial), then $P_j(x,y)$ is bounded by a scaled
convolution kernel (of the form on the right side of 
(\ref{pmx})) which
has rapid decay.  Spline wavelets (see \cite{St1, Ba1, 
Le1}) satisfy the conditions in part (a), and wavelets
constructed by P.~G.~Lemari\'e and Y.~Meyer \cite{Le-Me1} 
are 
an example of wavelets of rapid decay.


\section{Rates of convergence}

We will now give necessary and sufficient conditions on 
the basic
wavelet $\M^{\l}$ and on the scaling function $\F$, for 
given
supremum-norm rates of convergence of wavelet expansions.  
The
conditions on ${\M^{\l}}$ (here given in terms of 
membership in
homogeneous Sobolev spaces) can be translated into 
differentiability
and then moment conditions on $\M^{\l}$ (see the 
introduction).

\begin{theorem} Given a multiresolution analysis with 
either \rom{(i)} a
scaling function $\phi\in \RB$, \rom{(ii)} basic wavelets 
satisfying
$\psi^\lambda\ln (2+|x|)\in \RB$ for each $\lambda$, or   
\rom{(iii)} a kernel
for the basic projection $P$ satisfying $\abs{P(x,y)}\leq 
H(|x-y|)$ with $H\in  
\RB$, the following conditions \rom{((a)} to $(\mathrm 
e^\prime))$ are
equivalent\,\rom{:}

\vspace{.06in}

\rom{(a)} The multiresolution expansion \rom{(}see 
Definition \rom{1.2)}
yields pointwise order of approximation  $s$.

$(\mathrm a^\prime)$ The multiresolution expansion yields 
pointwise order of 
approximation  $s$ in every Sobolev space $H^r$ for $r 
\geq s+d/2$.

\rom{(b)} The projection $I-P_j\colon \ H_h^{s+
d/2}\rightarrow
L^\infty$ is a bounded operator, where $I$ is the identity 
and $d$
denotes dimension.

\vspace{.06in}

If there exists a family $\{ \psi^\lambda\}_\lambda \subset
\RB$ of basic wavelets corresponding to $\{ P_j\}$, 
then\,\rom{:}

\rom{(c)} For every family of basic wavelets
$\{\psi^\lambda\}_\lambda\subset \RB$ corresponding to $\{ 
P_j\}$,   
and for each $\lambda,\ \psi^\lambda\in H_h^{-s-d/2}$.

\rom{(d)} For every family of basic wavelets $\{
      \psi^\lambda\}_\lambda\subset \RB$  corresponding to 
$\{  
P_j\}$, and
       for each $\lambda$,
       \begin{equation}
       \int_{|\xi |<\epsilon} |\psi^\lambda (\xi )|^2 |\xi        
\label{intd}
        |^{-2(s+d/2)}\,d\xi <\infty 
      \end{equation}
     for some \rom{(}or for all\/\rom{)}  $\epsilon >0$.

$(\mathrm d^\prime)$  For some family of basic wavelets
$\{\psi^\lambda\}_\lambda \subset \RB$  corresponding to 
$\{P_j\}$,
      equation\ \rom{(}\ref{intd}\rom{)} holds.


If in addition there exists a scaling function $\phi\in\RB$
corresponding to some family of basic wavelets as above, 
then\,\rom{:}

\rom{(e)} For every such scaling function,
      \begin{equation}
      \int_{|\xi |<\epsilon}\left((2\pi)^d|\hat\phi (\xi 
)|^2 - 1\right)|\xi  
|^{-2(s+d/2)}\,d\xi  <\infty  \label{inte}
      \end{equation}
for some \rom{(}or all\/\rom{)} $\epsilon >0$.

$(\mathrm e^\prime)$ For some such scaling function 
$\phi$, equation
\rom{(}\ref{inte}\rom{)} holds.

\end{theorem}

The above conditions on the scaling function can also be 
given in
modified form in the case where the scaling function 
$\phi$ has
nonorthonormal integer translates.  Necessary and sufficient
conditions for convergence rates of spline and other 
nonorthonormal
expansions can then be obtained directly from the same 
arguments as
above. 

\vspace{.3in}

\section{Acknowledgments} 

The authors wish to thank Guido
Weiss for bringing us together. The first author would 
like to
acknowledge that this work was done under the direction of 
Richard
Rochberg and Mitchell Taibleson, and she thanks them for 
their
guidance. The last two authors wish to thank Ingrid 
Daubechies, Bjorn
Jawerth, Stephane Mallat, and Robert Strichartz for useful 
comments.




\begin{thebibliography}{AAAAAA}

\bibitem[Ba]{Ba1} G. Battle, {\it A block spin 
construction of ondelettes. 
Part \rom{I:} Lemari\'{e} functions}, {Comm. Math. Phys.} 
{\bf 110} (1987),
601--615.

\bibitem[BDR]{Bo-De-Ro} C. de Boor, R. DeVore, and A. Ron,
{\it Approximation from shift-invariant subspaces of 
$L^2({\bold R^d})$},
preprint.

\bibitem[BR]{Bo-Ro} C. de Boor and A. Ron, {\it Fourier 
analysis of the
approximation power of principal shift-invariant 
subspaces}, preprint.

\bibitem[Ca]{Ca1} L. Carleson, {\it On convergence and 
growth of partial sums of
Fourier series}, {Acta Math.} {\bf 116} (1966), 135--157.

\bibitem[Da1]{Da1} I. Daubechies, {\it Orthonormal bases 
of compactly 
supported wavelets}, {Comm. Pure Appl. Math.} {\bf 41}
(1988), 909--996. 

\bibitem[Da2]{Da2} \bysame,\ {\it Ten lectures on 
wavelets}, SIAM, 
Philadelphia, PA, 1992.

\bibitem[GM]{Gr-Mo1} A. Grossman and J. Morlet, {\it 
Decomposition of  
Hardy functions into square-integrable wavelets of 
constant shape,}  
SIAM J. Math. Anal. {\bf 15} (1984) 723--736.

\bibitem[Ha]{Ha1} A. Haar, {\it Zur Theorie der orthogonalen
Funktionensysteme}, {Math. Ann.} {\bf 69} (1910), 336.

\bibitem[Hu]{Hu1} R. A. Hunt, {\it On the convergence of 
Fourier series},
Proceedings of the Conference on Orthogonal Expansions and 
Their
Continuous Analogues (Southern Illinois University), 1968,
pp. 235--255.

\bibitem[Le]{Le1} P. G. Lemari\'e,  {\it Ondelettes \'a 
localisation
exponentielle}, J. Math. Pures Appl. (to appear).

\bibitem[LM]{Le-Me1}  P. G. Lemari\'e and Y. Meyer, {\it 
Ondelettes et  
bases Hilbertennes}, Rev. Math. Iberamericana {\bf 2} 
(1986), 1--18.

\bibitem[KKR]{Ke-Ko-Ra1} S. Kelly, M. Kon, and L. Raphael, 
{\it  
Pointwise convergence of wavelet expansions} (to appear).

\bibitem[Ma]{Ma1}  S. Mallat, {\it Multiresolution 
approximation and wavelets},
{Trans. Amer. Math. Soc.} {\bf 315} (1989), 69--88.

\bibitem[Me]{Me1} Y. Meyer,  {\it Ondelettes}, Hermann, 
Paris, 1990.

\bibitem[SF]{St-Fi} G. Strang and G. Fix, {\it A Fourier 
analysis of the finite
element variational method}, {Constructive Aspects of
Functional Analysis}, Edizioni Cremonese, Rome, 1973.  

\bibitem[St]{St1} J.-O. Str\"omberg,  {\it A modified 
Franklin system and
higher-order spline system on ${\bold R}^n$ as 
unconditional bases for
Hardy spaces}, Conference on Harmonic Analysis in Honor of 
Antoni
Zygmund, Vol. II (W. Beckner et al., eds.), Wadsworth Math. 
Ser., Wadsworth, Belmont, CA, 1981, pp. 475--494.

\bibitem[Wa1]{Wa1}G. Walter, {\it Approximation of the 
delta function by
wavelets}, preprint, 1992.

\bibitem[Wa2]{Wa2} \bysame,\ {\it Pointwise convergence of 
wavelet 
expansions}, preprint, 1992.

\end{thebibliography}
\end{document}